\documentclass[11pt]{article}
\usepackage{amssymb,amscd,amsthm,verbatim}

\usepackage{amsmath,amssymb,amsthm}
\usepackage{amsfonts}
\usepackage[mathscr]{eucal}
\begin{document}
\newtheorem{thm}{Theorem}
\numberwithin{thm}{section}
\newtheorem{lemma}[thm]{Lemma}
\newtheorem{remark}{Remark}
\newtheorem{corr}[thm]{Corollary}
\newtheorem{proposition}{Proposition}
\newtheorem{theorem}{Theorem}[section]
\newtheorem{deff}[thm]{Definition}
\newtheorem{case}[thm]{Case}
\newtheorem{prop}[thm]{Proposition}
\numberwithin{equation}{section}
\numberwithin{remark}{section}
\numberwithin{proposition}{section}
\newtheorem{corollary}{Corollary}[section]
\newtheorem{others}{Theorem}
\newtheorem{conjecture}{Conjecture}\newtheorem{definition}{Definition}[section]
\newtheorem{cl}{Claim}
\newtheorem{cor}{Corollary}
\newcommand{\ds}{\displaystyle}

\newcommand{\stk}[2]{\stackrel{#1}{#2}}
\newcommand{\dwn}[1]{{\scriptstyle #1}\downarrow}
\newcommand{\upa}[1]{{\scriptstyle #1}\uparrow}
\newcommand{\nea}[1]{{\scriptstyle #1}\nearrow}
\newcommand{\sea}[1]{\searrow {\scriptstyle #1}}
\newcommand{\csti}[3]{(#1+1) (#2)^{1/ (#1+1)} (#1)^{- #1
 / (#1+1)} (#3)^{ #1 / (#1 +1)}}
\newcommand{\RR}[1]{\mathbb{#1}}
\thispagestyle{empty}
\begin{titlepage}
\title{\bf    Higher order PDE's and iterated processes }
\author{Erkan Nane\thanks{ Supported in part by NSF Grant \# 9700585-DMS }\\
Department of Mathematics\\
Purdue University\\
West Lafayette, IN 47906 \\
enane@math.purdue.edu}
\maketitle
\begin{abstract}
\noindent {\it We  introduce a class of stochastic processes based
on symmetric $\alpha$-stable processes, for $\alpha \in (0,2]$.
 These are obtained by taking Markov processes and replacing the
 time parameter with the modulus of
 a symmetric $\alpha$-stable process. We call them $\alpha$-time
 processes. They generalize  Brownian time processes studied
in \cite{allouba1, allouba2, allouba3}, and they introduce new
interesting examples.  We establish the connection of
 $\alpha-$time processes to some higher order PDE's for $\alpha$ rational.
 We also study
 the exit problem for $\alpha$-time processes as they exit regular domains
  and connect them to
 elliptic PDE's. We also obtain the PDE connection of subordinate killed
 Brownian motion in bounded domains of regular boundary.
}
\end{abstract}
\end{titlepage}

\section{Introduction }

The link between concepts from probability and partial
differential equations (PDE's) helped solve problems in analysis
or find easier and shorter proofs for well-known results.
Researchers have been fascinated by these kinds of links. The
classical well-known connection of a PDE and a stochastic process
is the Brownian motion and heat equation connection. Let $X_{t}\in
\RR{R}^{n}$ be Brownian motion started at $x$. Then the function
$$
u(t,x)=E_{x}[f(X_{t})]
$$
solves the  Cauchy problem
\begin{eqnarray}
\frac{\partial}{\partial t}u(t,x)\ & = & \Delta u(t,x), \ \ \ \ \  \ t>0, \ \ x\in \RR{R}^{n}\nonumber\\
u(0,x)& = & f(x), \ \ \ \ \ x\in \RR{R}^{n}.\nonumber
\end{eqnarray}

In recent years, starting with the articles of Burdzy
\cite{burdzy1, burdzy2}, researchers had interest in iterated
processes in which one changes the time parameter with
one-dimensional Brownian motion, see \cite{allouba1, allouba3,
burdzy1, burdzy2, bukh, deblassie, funaki, koslew, xiao} and
references there in. The connections of iterated Brownian motion
(IBM) $Z_{t}=X(Y_{t})$, where $X$ is a two-sided Brownian motion
and $Y$ is another Brownian motion independent of $X$, and the
biLaplacian have been found by Allouba and Zheng \cite{allouba1}
in (2001) and by DeBlassie \cite{deblassie} in (2004). They showed
that the function
$$
u(t,x)=E_{x}[f(Z_{t})]
$$
solves the Cauchy problem
\begin{eqnarray}
\frac{\partial}{\partial t}u(t,x)\ & = &
\frac{{\Delta}f(x)}{\sqrt{\pi t}}\ + \
{\Delta}^{2}u(t,x),
\ \ \ \ \  \ t>0, \ \ x\in \RR{R}^{n}\nonumber\\
u(0,x)& = & f(x), \ \ \ \ \ x\in \RR{R}^{n}.\nonumber
\end{eqnarray}
The non-Markovian property of IBM is reflected by the appearance
of the initial function $f(x)$ in the PDE.

The results in the above mentioned articles are a motivation for
us to define and study our processes in this paper. Our aim in
this paper is two fold. In the first part we define $\alpha$-time
processes and study their PDE connection. In the second part, we
study the PDE connection of subordinate killed Brownian motion
over bounded domains of regular boundary.

In analogy with the Brownian time processes (BTP)'s studied in
\cite{allouba1, allouba2}, we define the `$\alpha$-time' processes
and establish the connection of these processes to several higher
order PDE's. Since the time clock that will replace ordinary time
is a symmetric $\alpha$-stable process, we first give the general
definition for these processes.

The $n$-dimensional symmetric stable process $X_{t}$ with index
$\alpha\in (0,2]$ is the process with stationary independent
increments whose transition density
$$
p_{t}^{\alpha}(x,y)=p^{\alpha}(t,x-y), \ \ \ \ \ (t,x,y)\in
(0,\infty)\times \RR{R}^{n}\times \RR{R}^{n},
$$
is characterized by the Fourier transform
$$
\int_{\RR{R}^{n}}
e^{iy.\xi}p^{\alpha}(t,y)dy=\exp(-t|\xi|^{\alpha}), \ \ \ \ \ t>0,
\xi \in \RR{R}^{n}.
$$
The process has right continuous paths, it is rotation and
translation invariant. For $\alpha =2$, this is Brownian motion.

 We now give the definitions of $\alpha$-time processes and some other iterated
processes. Let $\alpha\in (0,2]$. Let $Y(t)$ be a 1-dimensional
$\alpha$-stable process independent of the continuous Markov
process $X^{x}(t)$, started at 0, define the `$\alpha$-time'
process  $\alpha $TP by $Z^{x}_{Y}(t):=X^{x}(|Y(t)|)$. Similarly,
we define kE$\alpha$TP the excursions based  `$\alpha$-time'
process $Z^{x,k}_{Y,c}$; from $\alpha$TP's by breaking up the path
of $|Y(t)|$ into excursion intervals-maximal intervals $(r,s)$ of
time on which $|Y(t)|>0$-and, on each such interval, we pick an
independent copy of the Markov process $X^{x}$ from a finite or an
infinite collection. Let $X^{x,1}(t)\cdots X^{x,k}(t)$ be
independent copies of $X^{x}(t)$ starting from the point $x$. On
each excursion interval of $|Y(t)|$ use one of the $k$ copies
chosen at random. For $\alpha =2$, these are the Brownian time
processes defined in \cite{allouba1, allouba2, allouba3}.  When
the outer Markov process is a Brownian motion, $\alpha =2$ and
$k=2$, this is the famous iterated Brownian motion of Burdzy
\cite{burdzy1}. We will study PDE connection of these processes
extensively in sections 2-4.

In section 5, we will study subordinate killed Brownian motion,
which can be realized as a composition of Brownian motion killed on
the boundary of regular domain and a subordinator. Let
$D\subset\RR{R}^{n}$ be a bounded domain and let $X_{t}^{D}$ be
the Brownian motion killed upon exiting $D$. The subordinate
killed Brownian motion $Z_{t}^{D}$ is defined as the process
obtained by subordinating $X^{D}$ via the $\alpha /2$-subordinator
$T_{t}$. More precisely, let $Z_{\alpha}^{D}=X^{D}(T_{t})$, $t>0$.
Then $Z_{\alpha}^{D}$ is a Hunt process on $D$. If we use
$(P^{D}_{t})_{t\geq 0}$ to denote the semigroup of $X^{D}$ then
the semigroup $Q_{t}^{\alpha}$ of $Z_{\alpha}^{D}$ is given by
$$
Q_{t}^{\alpha}f(x)=\int_{0}^{\infty}P_{s}^{D}f(x)u^{\alpha
/2}_{t}(s)ds.
$$

Potential theory of this process has been studied by R. Song and
Z. Vodra\u{c}ek in \cite{songvod} over Lipschitz domains and
$C^{1,1}$ domains. In particular, they established intrinsic
ultraconractivity of the semigroup of this process and bounds on
the Green function and the jumping kernel.

R. Song in \cite{song}, established sharp bounds on the density,
 Green function and the jumping function of subordinate killed Brownian motion.
 For other properties of subordinate killed Brownian motion check the references in \cite{song} and \cite{songvod}.

The paper is organized as follows. We give the main results and
proofs of PDE connection of $\alpha$-time processes in \S2. We
study the exit problem from bounded domains in \S3. We collect
some useful lemmas in \S4. In \S5, we study the PDE connection of
subordinate killed Brownian motion in bounded domains of regular
boundary.

\section{ $\alpha$-time processes}

The first PDE connection of `$\alpha$-time' processes is the
 following. When $\alpha =1$, we call this process a Cauchy-time
 process (CTP) and denote it by $Z^{x}_{C}(t)$, for $C(t)$ a real-valued Cauchy process.

\begin{theorem}\label{cauchy}
Let ${\cal T}_{s}f(x)=E[f(X^{x}(s))]$ be the semigroup of the
continuous Markov process $X^{x}(t) $ and let ${\cal A}$ be its
generator. Let $\alpha=1$. Let $f$ be a bounded measurable
function in the domain of ${\cal A}$, with $D_{ij}f$ bounded and
H\"{o}lder continuous for all $1\leq i,\ j \leq n$. If
$u(t,x)=E[f(Z^{x,k}_{C,c}(t))]$ for any $k\in \RR{N}$, then $u$
solves
\begin{eqnarray}
\frac{\partial ^{2}}{\partial t^{2}}u(t,x)\ &\ = \ & -\frac{2{\cal A} f(x)}{\pi t}\ - \ {\cal A}^{2}u(t,x),
\  \  \ t>0, \ \ x\in \RR{R}^{n}\label{PDE-CONNECT2}\\
u(0,x) &\  =\  & \ f(x), \ \ \ \ \ x \in \RR{R}^{n},\nonumber
\end{eqnarray}
where the operator ${\cal A} $ acts on $u(t,x)$ as a function of
$x$ with $t$ fixed. In particular, if $X^{x}(t) $ is Brownian
motion started at $x$ and $\Delta$ is the standard Laplacian, then
$u$ solves
\begin{eqnarray}
\frac{\partial^{2}}{\partial t^{2}}u(t,x)\ & = &-\frac{2\Delta
f(x)}{\pi t}\ - \ \Delta ^{2} u(t,x),
\ \ \ \ \  \ t>0, \ \ x\in \RR{R}^{n}\label{PDE-CONNECT3}\\
 u(0,x) & = & f(x), \ \ \ \ \ x\in \RR{R}^{n}.\nonumber
\end{eqnarray}
\end{theorem}

\begin{remark}
Notice the effect of the initial function is different than the
Brownian motion heat equation case, it is very similar to the
effect of the initial function in IBM-PDE connection. The
appearance of the initial function in the PDE reflects the
non-Markovian property of the iterated processes. This PDE is very
similar to the one in \cite[Theorem 0.1]{allouba1} in that this
has the same order space derivatives (i.e Laplacian and
biLaplacian), but on the other hand this has second order time
derivative and  the coefficients are different.
\end{remark}

\begin{proof}
 We first consider the process $Z^{x}_{C}$. We use
the representation
$$
u(t,x)=E[f(Z^{x}_{C}(t))] =2\int_{0}^{\infty}p_{t}^{1}(0,s){\cal
T}_{s}f(x)ds,$$ where $p_{t}^{1}(0,s)=\frac{t}{\pi(s^{2}+t^{2})}$
is the transition density of the Cauchy process on $\RR{R}$. Using
dominated convergence and the fact that $p_{t}^{1}(0,s)$ satisfy
$$
(\frac{\partial^{2}}{\partial s^{2}}+\frac{\partial^{2}}{\partial
t^{2}})p_{t}^{1}(0,s)=0,
$$
we obtain
\begin{eqnarray}
\frac{\partial^{2}}{\partial t^{2}}E[f(Z^{x}_{C}(t))]& = &2\int_{0}^{\infty}\frac{\partial^{2}}{\partial t^{2}}p_{t}^{1}(0,s){\cal T}_{s}f(x)ds\nonumber\\
&= &2  \int_{0}^{\infty}-\frac{\partial^{2}}{\partial
s^{2}}p_{t}^{1}(0,s){\cal T}_{s}f(x)ds.
\end{eqnarray}

We now integrate by parts twice and observe that boundary terms
always vanish at $\infty$ (as $s\to \infty$) and we have
$\frac{\partial}{\partial s}p_{t}^{1}(0,0)=0$, but
$p_{t}^{1}(0,0)>0$. Thus
\begin{eqnarray}
\frac{\partial^{2}}{\partial t^{2}}E[f(Z^{x}_{C}(t))]& = &2\int_{0}^{\infty}\frac{\partial}{\partial s}p_{t}^{1}(0,s)\frac{\partial}{\partial s}{\cal T}_{s}f(x)ds\nonumber\\
&=&-2p_{t}^{1}(0,0){\cal A}f(x)
-2\int_{0}^{\infty}p_{t}^{1}(0,s){\cal A}^{2}{\cal
T}_{s}f(x)ds.\nonumber
\end{eqnarray}
Taking the application of ${\cal A}^{2}$ outside the integral using the conditions on $f$ and $D_{ij}f$ by Lemma \ref{cauchy-diff-int}
  we get
$$
\frac{\partial ^{2}}{\partial t^{2}}u(t,x) = -2p_{t}^{1}(0,0){\cal
A} f(x)\ - \ {\cal A}^{2}u(t,x).
$$

For the other processes we have by similar arguments as in Theorem
0.1 in Allouba and Zheng \cite{allouba1} considering the excursion
intervals of $|C(t)|$
$$
E[f(Z^{x}_{C}(t))]=E[f(Z^{x,k}_{C,c}(t))].
$$
We provide the argument in Allouba and Zheng for completeness. Let $e^{-}(t)$ be the
 $|C(t)|$-excursion immediately preceding the
excursion straddling $t$, $e(t)$; condition on the event that we
pick the $j$th copy of $X^{x}$ on $e^{-}(t)$ (uniformly among the
$k$ copies of $X^{x}$), using the independence of the process
$X^{x,j}$ on $e^{-}(t)$ from $C(t)$ and from the following choice
of $X^{x}$ copy on $e(t)$ to get
\begin{eqnarray}
E[f(Z^{x,k}_{C,c}(t))]&= &
2\sum_{j=1}^{k}\int_{0}^{\infty}p_{t}^{1}(0,s){\cal
T}_{s}f(x)P[\text{picking jth copy on}
\  e^{-}(t) ] ds\nonumber\\
& = &  \frac{2}{k}\sum_{j=1}^{k}\int_{0}^{\infty}p_{t}^{1}(0,s){\cal T}_{s}f(x) ds =
 2\int_{0}^{\infty}p_{t}^{1}(0,s){\cal T}_{s}f(x) ds\nonumber\\
& = & E[f(Z^{x}_{C}(t))].
\end{eqnarray}
This concludes the proof of
Theorem \ref{cauchy}.
\end{proof}

Next, we solve a similar PDE as in \cite[Theorem 1.2]{allouba2}
which is obtained by running an $\epsilon$-scaled CTP and
averaging the product of $f(Z_{\epsilon C,c}^{x,k}(t))$ with the
negative exponential of $|C(t)| / \epsilon$, for $C(t)$ a
real-valued Cauchy process. This looks like the Feynman-Kac
formula when $\epsilon =1$. We state this theorem since it is a
step towards the probabilistic study of Linearized Cahn-Hilliard
and Kuramoto-Sivashinsky type PDE's as in \cite{allouba2}.

\begin{theorem}\label{theorem3}
Under the same conditions on $f$ as in Theorem \ref{cauchy}, and
for $\epsilon >0$, If
\begin{equation}\label{theorem3-1}
u_{\epsilon}(t,x)=E\left[ f(Z_{\epsilon C,c}^{x,k}(t))\exp\left(
 -\frac{|C(t)|}{\epsilon}\right) \right]
\end{equation}
for any $k\in \RR{N}$, then
$u_{\epsilon}$ solves
$$
\left\{
\begin{array}{rll}
\frac{\partial^{2}}{\partial t^{2}}u_{\epsilon}(t,x) & =
-\frac{2}{\pi t}\left[ \epsilon {\cal A}
f(x)-\frac{1}{\epsilon}f(x) \right] -
\frac{1}{\epsilon^{2}}u_{\epsilon} (t,x) & \\
& & \\
 & +  2{\cal
A}u_{\epsilon}(t,x)- \epsilon^{2}{\cal
A}^{2}u_{\epsilon}(t,x), & t>0,  \ x\in \RR{R}^{n},  \\
& & \\
 u_{\epsilon}(0,x) &=  f(x)=\lim_{t\downarrow 0,\ y\to x}
u_{\epsilon}(t,y), & x\in   \RR{R}^{n}. \label{theorem3-2}
\end{array}
\right.
$$

In particular, If the outer process $X^{x}$ in
(\ref{theorem3-1}) is a Brownian motion, then $u_{\epsilon}(t,x)$
solves
$$
\left\{
\begin{array}{rll}
\frac{\partial^{2}}{\partial t^{2}}u_{\epsilon}(t,x) & =
-\frac{2}{\pi t}\left[ \epsilon \Delta
f(x)-\frac{1}{\epsilon}f(x) \right] -
\frac{1}{\epsilon^{2}}u_{\epsilon} (t,x) & \\
& & \\
 & +  2 \Delta u_{\epsilon}(t,x)- \epsilon^{2}\Delta^{2}u_{\epsilon}(t,x), & t>0,  \ x\in \RR{R}^{n},  \\
& & \\
 u_{\epsilon}(0,x) &=  f(x)=\lim_{t\downarrow 0,\ y\to x}
u_{\epsilon}(t,y), & x\in   \RR{R}^{n}.
\end{array}
\right.
$$
\end{theorem}

\begin{remark}
This time the initial function affects the PDE differently(i.e
through both $f$ and ${\cal A}f$). We also comment that, the last
two terms in the PDE (the biLaplacian and the Laplacian of the
solution $u_{\epsilon}$) look like those in a linearized
Cahn-Hilliard equation with the correct $\epsilon$-scaling, but
with the opposite sign for $\Delta$.
\end{remark}

\begin{proof}
It is enough to prove the CTP case. Let
\begin{equation}\label{theorem3-4}
u_{\epsilon}(t,x)=E\left[ f(Z_{\epsilon C}^{x}(t))\exp\left(
-\frac{|C(t)|}{\epsilon}\right) \right]
\end{equation}
and
\begin{equation}\label{theorem3-5}
v_{\epsilon}(s,x)=E\left[ f(X^{x}(\epsilon s))\exp\left(
-\frac{s}{\epsilon}\right)\right]=\exp\left(
-\frac{s}{\epsilon}\right){\cal T}_{\epsilon s}f(x).
\end{equation}

We then have
\begin{equation}\label{theorem3-6}
u_{\epsilon}(t,x)=2\int_{0}^{\infty} v_{\epsilon}(s,x)
p_{t}^{1}(0,s)ds.
\end{equation}
The rest of the proof is similar to the proof of Thereom 1.2 in
Allouba \cite{allouba2}.
\end{proof}

The next result gives a Feynman-Kac type formula for CTP's and
connects it to fourth order PDE's.

\begin{theorem}
Assume that $f, c: \RR{R}^{n}\to \RR{R}$ are bounded, $c\leq 0$,
and $D_{ij}f$ and $D_{ij}c$ are bounded and H\"{o}lder continuous
with exponent $0<\beta \leq 1$, for $1\leq i,j\leq n$. If the
$|D_{i,j}v(s,x)|\leq K_{T}$, for all $(s,x)\in [0,T]\times
\RR{R}^{n}$, for any time $T>0$, for all $i,j$, where $K_{T}>0$ is
a constant depending only on $T$ and
\begin{equation}\label{theorem4-1}
v(s,x)=E\left[ f(X^{x}(s))\exp\left( \int_{0}^{s}
c(X^{x}(r))dr\right)  \right]
\end{equation}
where $X^{x}$ is an $n-$dimensional Brownian motion starting
at $x$, then
\begin{equation}\label{theorem4-2}
u(t,x)=E \left[ f(Z_{C}^{x}(t))\exp \left(
\int_{0}^{|C(t)|}c(X^{x}(r))dr \right) \right]
\end{equation}
solves
\begin{equation}
 \left\{
\begin{array}{rl}
\frac{\partial^{2}}{\partial t^{2}}u(t,x) =& -\frac{2}{\pi
t}\left[\Delta f(x)+ c(x)f(x)\right]  \\
 & \\
&   -\left[ \Delta c(x)+ c^{2}(x)\right]u(t,x) -2\nabla
c(x).\nabla
u(t,x)  \\
&  \\
 & -  c(x)\Delta u(t,x) -\Delta^{2}u(t,x), \
  \ \ \ \ t>0, \ x\in \RR{R}^{n},\\
  &\\
    u(0,x)= & f(x)\  = \ \lim_{t\downarrow 0,\ y\to x} u(t,y),  \ \  \ \ \ \ \ \
x\in \RR{R}^{n}.
\end{array}
\right.
\end{equation}
\end{theorem}

\begin{remark}
Again the effect of the initial function is through $\Delta$ and
$\Delta ^{2}$. The effect of the initial function fades away as
$t\to\infty$ at a rate of $2 / \pi t$.
\end{remark}

\begin{proof}
Let $u$ and $v$ be defined as in (\ref{theorem4-2}) and
(\ref{theorem4-1}), respectively.
 Then
 \begin{equation}\label{theorem4-3}
 u(t,x)=2\int_{0}^{\infty} p_{t}^{1}(0,s)v(s,x) ds.
 \end{equation}
 Since $X^{x}$ is Brownian motion starting at $x$, we have
 \begin{equation}\label{tehorem4-4}
\frac{\partial}{\partial s} v(s,x)=\Delta v(s,x)+
c(x)v(s,x),\  \mathrm{in }\ (0,\infty)\times \RR{R}^{n}.
 \end{equation}
The rest of the proof is similar to the proof of Thereom 1.3 in
Allouba \cite{allouba2}.
\end{proof}

Now  we define an iterated process motivated by the article of
Allouba \cite{allouba3} that leads to the study of  Linearized
Kuramoto-Sivashinsky PDE. We call the process as
imaginary-Cauchy-time-Brownian-angle process (ICTBAP). The
definition needs the introduction of complex imaginary number
$i=\sqrt{-1}$.

Let $f:\RR{R}^{n}\to\RR{R}$
$$
{\cal B} _{C}^{f,X}(t,x)=\left\{ \begin{array}{ll}
f(X^{x}(iC(t)))\exp(iC(t)), &C(t)\geq 0,\\
f(X^{-ix}(-iC(t)))\exp(iC(t)), &C(t)<0.
\end{array}
\right.
$$
where  $X^{x}$ is an $\RR{R}^{n}$-valued Brownian motion starting
from $x\in \RR{R}^{n}$, $X^{-ix}$ is an independent
$i\RR{R}^{n}$-valued BM starting at $-ix$, and both are
independent of the inner $\RR{R}$-valued Cauchy process $C$. We
think of the imaginary-time processes $\{ X^{x}(is),\ s\geq 0 \}$
and $\{ iX^{-ix}(-is),\ s\leq 0 \}$ as having the same complex
Gaussian distribution on $\RR{R}^{n}$ with the corresponding
complex distributional density
$$
p_{is}^{n}(x,y)=\frac{1}{(4\pi is)^{n/2}}e^{-|x-y|^{2}/4is} .
$$

The space derivative terms in the PDE in the following Theorem is
the same as in a linearized Kuramoto-Sivashinsky PDE.

\begin{theorem}
Let $f\in C_{c}^{2}(\RR{R}^{n};\RR{R})$ with $D_{ij}f$ H\"{o}lder
continuous with exponent $0<\beta \leq 1$, for all $1\leq i,j \leq
n$. Let
$$
v(s,x)= \exp(is) \int_{\RR{R}^{n}} f(y) \frac{1}{(4\pi is)^{n/2}
}e^{-|x-y|/4is}dy,
$$
$$
u(t,x)=\int_{-\infty}^{0}v(s,x)p_{t}^{1}(0,s)ds
+\int_{0}^{\infty}v(s,x)p_{t}^{1}(0,s)ds,
$$
where
$$
p_{t}^{1}(0,s)=\frac{t}{\pi (s^{2}+t^{2})}
$$
is the transition density of the inner (one-dimensional) Cauchy
process. Then $u(t,x)$ solves 
$$
\left\{\begin{array}{ll} \frac{\partial ^{2}}{\partial
t^{2}}u(t,x) = \Delta ^{2}u(t,x)+2\Delta u(t,x)+u(t,x),
 & t>0, x\in \RR{R}^{n}, \\
 & \\
u(0,x) =  f(x),  & x\in \RR{R}^{n}.
\end{array}
\right.
$$
\end{theorem}

\begin{proof}
We comment that
$$
v(s,x)=E[f(X^{x}(is))\exp(is)]
$$
and
$$
u(t,x)=E[{\cal B} _{C}^{f,X}(t,x)].
$$

The proof is the same  as in \cite{allouba3} except the start point that $p_{t}^{1}(0,s)$
satisfy

$$
(\frac{\partial^{2}}{\partial s^{2}}+\frac{\partial^{2}}{\partial
t^{2}})p_{t}^{1}(0,s)=0,
$$
which explains the sign change from Theorem 1.1 in
\cite{allouba3}. We can differentiate under the integral by
Lemma 2.1 in \cite{allouba3} and Lemma \ref{cauchy-diff-int}. The
rest of the proof is similar to the proof of Thereom 1.1 in
Allouba \cite{allouba3}.
\end{proof}

For  $\alpha \neq 1$, the PDE is more complicated since  kernels
of symmetric $\alpha$-stable processes satisfy a higher order PDE.
We also have to assume that we can integrate under the integral as
much as we need in the case where the outer process is BM (or in
general we can take the operator out of the integral). This is
valid for $\alpha=1/m$, $m=2,3,\cdots$ by Lemma
\ref{alpha-interchange} below.

\begin{theorem}\label{alphapde}
Let $\alpha \in (0,2)$ be rational $\alpha=l/m$, where $l$ and $m$
are relatively prime. Let ${\cal T}_{s}f(x)=E[f(X^{x}(s))]$ be the
semigroup of the continuous Markov process $X^{x}(t) $ and let
${\cal A}$ be its generator. Let $f$ be a bounded measurable
function in the domain of ${\cal A}$, with $D^{\gamma}f$  bounded and
H\"{o}lder continuous for all multi index $\gamma$ such that $|\gamma |=l$. If
$u(t,x)=E[f(Z^{x,k}_{Y,c}(t))]$ for any $k\in \RR{N}$, then $u$
solves the PDE
\begin{eqnarray}
(-1)^{l+1}\frac{\partial ^{2m}}{\partial t^{2m}}u(t,x) & =  &
-2\sum_{i=1}^{l}\left(
\frac{\partial^{2l-2i}}{\partial s^{2l-2i}}p_{t}^{\alpha}(0,s)|_{s=0}\right){\cal A}^{2i-1}f(x)\   \nonumber\\
&\ &\ -\ {\cal A}^{2l}u(t,x),\ \ \ t>0,\ \  x\in \RR{R}^{n} \nonumber\\
u(0,x) &\  =\  & \ f(x), \ \ \ \ \ x \in \RR{R}^{n},\nonumber\
\end{eqnarray}
where the operator ${\cal A} $ acts on $u(t,x)$ as a function of
$x$ with $t$ fixed. In particular, if $X^{x}(t) $ is Brownian
motion started at $x$ and $\Delta$ is the standard Laplacian, then
$u$ solves
\begin{eqnarray}
(-1)^{l+1}\frac{\partial ^{2m}}{\partial t^{2m}}u(t,x)\ &\ = \ &
-2\sum_{i=1}^{l} \left(\frac{\partial^{2l-2i}}{\partial
s^{2l-2i}}p_{t}^{\alpha}(0,s)|_{s=0}\right) \Delta^{2i-1}f(x)\ 
  \nonumber\\
 &\ &\ -\ \Delta^{2l} u(t,x), \ \ \ \ t>0,\ \   x\in \RR{R}^{n} \nonumber\\
u(0,x) &\  =\  & \ f(x), \ \ \ \ \ x \in \RR{R}^{n}.\nonumber\
\end{eqnarray}
\end{theorem}

\begin{proof}
We use the representation
$$
u(t,x)=E[f(Z^{x}_{Y}(t))]=2\int_{0}^{\infty}p_{t}^{\alpha}(0,s){\cal
T}_{s}f(x)ds
$$
and the fact that the transition density of the process $Y(t)$
satisfies
$$
(\frac{\partial^{2}}{\partial
s^{2}})^{l}+(-1)^{l+1}\frac{\partial^{2m}}{\partial
t^{2m}})p_{t}^{\alpha}(0,s)=0,\ \ \ \ \ (t,x)\in (0,\infty)\times
\RR{R}^{n},
$$
 from Lemma \ref{alphalemma}. Then, we use integration by parts as many as
 needed to pass the $s$ derivatives to ${\cal T}_{s}f(x)$ from $p_{t}^{\alpha}(0,s)$.
 The boundary terms are all zero at $\infty$ (as $s$ goes to $\infty$). And also all the odd $s$ derivatives of $p_{t}^{\alpha}(0,s)$ at $0$ are zero.
\end{proof}

\section{Exit problem from bounded domains}

In this section, we study the exit problem for the Cauchy-time processes. Let
$G\subset \RR{R}^{n}$ be a bounded open set with regular boundary.
Let $T_{G}^{x}=\inf \{t\geq 0; \ Z^{x}_{Y}(t)\notin G\}$  be the
first exit time of CTP from $G$. 





The CTBM's differ from the BTBM's (i.e. the outer process is
Brownian motion) in the following Theorem. The BTBM solves a
fourth order PDE, but the CTBM solves the usual second order
elliptic PDE.

\begin{theorem}\label{elliptic2}
Let $D\subset \RR{R}^{n}$ be a domain with a regular boundary. Let
the outer process in $Z_{C}^{x}$ be a Brownian motion started at
$x$ and let $T_{D}^{x}$ be the first exit time of $Z_{C}^{x}$ from
the interior of $D$. If $u(x)=E_{x}[T_{D}^{x}]$, then $u$
satisfies
\begin{eqnarray}
\Delta u(x)\ & = & -1 \ \ \  x\in D\nonumber\\
 u(x)& = & 0, \ \ \ \ \ x\in \partial G.\nonumber
\end{eqnarray}
\end{theorem}

\begin{proof}
 Here the outer process in $Z_{Y}^{x}$ is a Brownian motion
started at $x$. Let $u(x)=E[T_{G}^{x}]$ and observe that
\begin{eqnarray}
T_{G}^{x}=\inf\{t\geq 0;\ \  Z_{Y}^{x}\notin  G\}&=& \inf \{t\geq 0;\ \ |Y(t)|\notin [ 0, \tau_{G}^{x}) \}\nonumber\\
\ &=& \inf \{ t\geq 0;\ \ Y(t)\notin
(-\tau_{G}^{x},\tau_{G}^{x})\},\nonumber
\end{eqnarray}
where $\tau_{G}^{x}=\inf\{ t\geq  0;\ \ X^{x}\notin G\}$. Thus
conditioning on $\tau_{G}^{x}$ we easily get
\begin{equation}
u(x)=E[E[T_{G}^{x}|\tau_{G}^{x}]]=E[\tau_{G}^{x}],
\end{equation}
from Getoor \cite{getoor}. But it is well-known that
$$u(x)=E[\tau_{G}^{x}]$$ satisfies
$$
\Delta u(x)=-1,
$$
for any bounded regular domain and $u(x)=0$ for $x\in \partial G$.
\end{proof}

\section{Useful Lemmas}

We start with a result that is useful in taking the derivative
under the integral. With the following lemma we  can now
establish  similar results to the ones in \cite{allouba1} and
\cite{allouba2}.

\begin{lemma}\label{cauchy-diff-int}
Let $X^{x}$ be an $n-$dimensional Brownian motion starting at $x$,
and let $f,g:\RR{R}^{n}\to \RR{R}$ be bounded and measurable such
that $D_{ij}f$ and $D_{ij}g$ are H\"{o}lder continuous, with
exponent $0<\beta \leq 1$, for $1\leq i,j\leq n$. Let
$$
u_{1}(t,x)=\int_{0}^{\infty}E[f(X^{x}(s))]p_{t}^{1}(0,s)ds,
$$
$$
u_{2}(t,x)=\int_{0}^{t}\int_{0}^{\infty}E[g(X^{x}(s))]p_{t}^{1}(0,s)dsdr.
$$
Then $\Delta ^{2}u_{1}(t,x) $ and $\Delta ^{2}u_{2}(t,x) $ are
finite and
$$
\Delta ^{2}u_{1}(t,x)=\int_{0}^{\infty}\Delta
^{2}E[f(X^{x}(s))]p_{t}^{1}(0,s)ds,
$$
$$
\Delta ^{2}u_{2}(t,x)=\int_{0}^{t}\int_{0}^{\infty}\Delta
^{2}E[g(X^{x}(s))]p_{t}^{1}(0,s)dsdr.
$$
\end{lemma}

\begin{proof}
The proof follows from the proof of Lemma 2.1 in \cite{allouba2}
with only changing the density of time process
$$
p_{t}^{1}(0,s)=\frac{t}{\pi(s^{2}+t^{2})},
$$
which is the transition density of the real-valued Cauchy process.
Our results follow from the facts,
$$
\int_{0}^{\infty}p_{t}^{1}(0,s)s^{1-\beta
/2}ds=\int_{0}^{\infty}\frac{t}{\pi(s^{2}+t^{2})s^{1-\beta /2}}ds
<\infty,
$$
and
$$
\int_{0}^{t}\int_{0}^{\infty}p_{r}^{1}(0,s)/s^{1-\beta
/2}dsdr=\int_{0}^{t}\int_{0}^{\infty}\frac{r}{\pi(s^{2}+r^{2})s^{1-\beta
/2}}dsdr <\infty.
$$
\end{proof}

For $\alpha\in (0,2)$ and $\alpha\neq 1$, we have the following
Lemma from DeBlassie \cite[Theorem 1.1]{deblassie2}.

\begin{lemma}\label{alphalemma}
Let $\alpha=l/m$ with $l,m$ relatively prime. The transition
density $p_{t}^{\alpha}(x,y)$ of the $\alpha$-stable process
satisfies the PDE for $y$ fixed
$$
(\Delta^{l}+(-1)^{l+1}\frac{\partial^{2m}}{\partial
t^{2m}})p_{t}^{\alpha}(x,y)=0,\ \ \ \ \ (t,x)\in (0,\infty)\times
\RR{R}^{n}.
$$
\end{lemma}

The next lemma allows us to take the biLaplacian out of the integral in Theorem \ref{alphapde} for some
values of the index of the symmetric stable process.
\begin{lemma}\label{alpha-interchange}
Let $X^{x}$ be an $n-$dimensional Brownian motion starting at $x$,
and let $f,g:\RR{R}^{n}\to \RR{R}$ be bounded and measurable such
that $D_{ij}f$ and $D_{ij}g$ are H\"{o}lder continuous, with
exponent $0<\beta \leq 1$, for $1\leq i,j\leq n$. Let $\alpha =1/m$, $m=2,3,\cdots$ and
$$
u(t,x)=\int_{0}^{\infty}E[f(X^{x}(s))]p_{t}^{\alpha}(0,s)ds.
$$
Then $\Delta ^{2}u(t,x) $ is
finite and
$$
\Delta ^{2}u(t,x)=\int_{0}^{\infty}\Delta
^{2}E[f(X^{x}(s))]p_{t}^{\alpha}(0,s)ds.
$$
\end{lemma}

\begin{proof}
The proof follows similarly to the proof of Lemma
\ref{cauchy-diff-int}, by showing that
$$
\int_{0}^{\infty}p_{t}^{\alpha}(0,s)s^{1-\beta
/2}ds=\int_{0}^{\infty}\int_{0}^{\infty} \frac{p_{l}(0,s) u^{\alpha /2}_{t}(l) }{ s^{1-\beta /2}}dlds
<\infty.
$$
Above we use the fact that
$$
p_{t}^{\alpha}(0,s)= \int_{0}^{\infty}p_{l}(0,s) u^{\alpha /2}_{t}(l) dl,
$$
where $p_{l}(0,s)$ is the density of real-valued Brownian motion,
and  $u^{\alpha /2}_{t}(l)$ is the density of the $\alpha /2$
subordinator (see \cite{deblassie2}  for this representation and
the bounds for integrals of $u^{\alpha /2}_{t}(l)$).
\end{proof}

\section{Subordinate killed Brownian motion}

Let $D\subset\RR{R}^{n}$ be a bounded domain and let $X_{t}^{D}$
be the Brownian motion killed upon exiting $D$. The subordinate
killed Brownian motion $Z_{t}^{D}$ is defined as the process
obtained by subordinating $X^{D}$ via the $\alpha /2$-subordinator
$T_{t}$. More precisely, let $Z_{\alpha}^{D}=X^{D}(T_{t})$, $t>0$.
Then $Z_{\alpha}^{D}$ is a Hunt process on $D$. If we use
$(P^{D}_{t})_{t\geq 0}$ denote the semigroup of $X^{D}$ then the
semigroup $Q_{t}^{\alpha}$ of $Z_{\alpha}^{D}$ is given by
$$
Q_{t}^{\alpha}f(x)=\int_{0}^{\infty}P_{s}^{D}f(x)u^{\alpha
/2}_{t}(s)ds.
$$

The following theorem establishes the PDE connection of subordinate killed Brownian motion.

\begin{theorem}\label{subordinate}
Let $0< \alpha =\frac{k}{m} <2$ where $k, m$ are relatively prime
integers. Let $f:D\to \RR{R}$ be a bounded function, then
$u(t,x)=Q_{t}^{\alpha}f(x)$ solves the following PDE,
$$
\left\{\begin{array}{ll} \Delta^{k} u(t,x)+(-1)^{k+1}\frac{\partial^{2m} }{\partial
t^{2m}}u(t,x)=0, & (t,x)\in (0,\infty)\times D \\
& \\
u(0,x) =  f(x),  & x\in D\\
& \\
 u(t,x)=0, & x\in \partial D.
\end{array}
\right.$$
\end{theorem}

To prove Theorem \ref{subordinate}, we need a differentiating
under the integral lemmas.
\begin{lemma}\label{laplacian-int}
Let $u(t,x)$ be as in Theorem \ref{subordinate}, then
$$
\Delta ^{k} u(t,x)=\int_{0}^{\infty}\Delta_{x}
^{k}P_{s}^{D}f(x)u^{\alpha /2}_{t}(s)ds.
$$
\end{lemma}

\begin{proof}
We know that
$$
\Delta ^{k}P^{D}_{s}f(x)=\sum_{l=1}^{\infty}(-\lambda_{l})^{k}
e^{-s \lambda _{l}} \varphi_{l}(x)\int_{D}\varphi_{l}(y)dy.
$$
Now
\begin{eqnarray}
\Delta ^{k} u(t,x)& = & \int_{0}^{\infty}\Delta
^{k}P_{s}^{D}f(x)u^{\alpha /2}_{t}(s)ds\nonumber \\
&=& \sum_{l=1}^{\infty}(-\lambda_{l})^{k}
\varphi_{l}(x)\int_{D}\varphi_{l}(y)dy \int_{0}^{\infty}e^{-s
\lambda _{l} }u^{\alpha /2}_{t}(s)ds\nonumber\\
& = &
\sum_{l=1}^{\infty}(-\lambda_{l})^{k}e^{-t(\lambda_{l})^{\alpha/2}}
\varphi_{l}(x)\int_{D}\varphi_{l}(y)dy, \label{laplace-trans}
\end{eqnarray}
where (\ref{laplace-trans}) follows because the Laplace transform
of the density of $T_{t}$ is
$$
\int_{0}^{\infty}e^{-s\lambda }u^{\alpha
/2}_{t}(s)ds=e^{-t\lambda^{\alpha /2}}, \ \ for \ \ \lambda >0.
$$
This last series is absolutely and uniformly convergent. Indeed
let $\delta >0$; since $p_{s}^{D}(x,y)\leq (4\pi s)^{-n/2}$, for
all $s>0$ we see that
\begin{eqnarray}
|e^{-s\lambda_{l}}\varphi_{l}(x)|&=&
|\int_{D}p_{s}^{D}(x,y)\varphi_{l}(y)dy|\nonumber\\
&\leq  &(4\pi s)^{-n/2}|| \varphi_{l}||_{2}(vol(D))^{1/2}\nonumber\\
& =& (4\pi s)^{-n/2}(vol(D))^{1/2},\nonumber
\end{eqnarray}
taking $s=1/\lambda_{l}$ gives that
$$
|\varphi_{l}(x)|\leq
e(4\pi)^{-n/2}\lambda_{l}^{n/2}[vol(D)]^{1/2},
$$
for all $x\in D$. Now since the volume of the domain $D $ is
finite, there exists $C>0$ such that
$$
\sum_{l=1}^{\infty}(-\lambda_{l})^{k}e^{-t(\lambda_{l})^{\alpha/2}}
\varphi_{l}(x)\int_{D}\varphi_{l}(y)dy\leq
C\sum_{l=1}^{\infty}e^{-\delta \lambda_{l}^{\alpha/2}/2},
$$
for all $(t,x)\in[ 0,\infty)\times D$ and the last series is
convergent by the Weyl's asymptotic lemma (see \cite{chavel}) we
have $c_{n,D}l^{n/2}\leq \lambda_{l}$ where $c_{n,D} $ is a
constant that depends on $n$ and the volume of the domain $D$.
\end{proof}

\begin{lemma}\label{t-derivative}
 Suppose $f$ is bounded. Set $u(t,x)=Q_{t}^{\alpha}f(x)$. Then for
 any integer $q\geq 0 $,
 \begin{equation}
 \frac{\partial^{q}}{\partial
 t^{q}}u(t,x)=\int_{D}\int_{0}^{\infty}f(y)p_{s}^{D}(x,y)\frac{\partial^{q}}{\partial
 t^{q}}u^{\alpha /2}_{t}(s)dsdy.
 \end{equation}
 \end{lemma}

 \begin{proof}
 Proof is the same as in Lemma 2.1 in \cite{deblassie2} almost
 word for word except we use the fact that
 $\int_{D}p_{s}^{D}(x,y)dy=1$.
 \end{proof}

 The next lemma is Lemma 3.1 in \cite{deblassie2}, we give it for
 completeness.

\begin{lemma}\label{density-pde}
 Let $\alpha=k/m$. For $s$ and $t$,
 $$
\left( \frac{\partial ^{k}}{\partial s^{k}} - \frac{\partial ^{2m}}{\partial
t^{2m}}\right)u^{\alpha /2}_{t}(s)=0.
 $$
\end{lemma}

\begin{proof}[Proof of Theorem \ref{subordinate}]
We first apply Lemma \ref{laplacian-int}, then integration by
parts repeatedly, and appealing equations (2.7)-(2.10) in
\cite{deblassie2}, to see the boundary terms are all $0$.
\begin{eqnarray}
\Delta ^{k} u(t,x)& = & \int_{D}\int_{0}^{\infty}f(y)\Delta_{x}
^{k}p_{s}^{D}(x,y)u^{\alpha /2}_{t}(s)dsdy \nonumber\\
& = &\int_{D}\int_{0}^{\infty}f(y) \left[
\frac{\partial^{k}}{\partial s^{k}}p_{s}^{D}(x,y) \right]
u^{\alpha /2}_{t}(s)dsdy\nonumber\\
& = &\int_{D}-f(y)\int_{0}^{\infty} \left[
\frac{\partial^{k-1}}{\partial s^{k-1}}p_{s}^{D}(x,y) \right]
\frac{\partial}{\partial s} u^{\alpha /2}_{t}(s)dsdy\nonumber\\
& =&\cdots =(-1)^{k}\int_{D}\int_{0}^{\infty}f(y)
p_{s}^{D}(x,y)\frac{\partial^{k}}{\partial s^{k}}u^{\alpha
/2}_{t}(s)dsdy. \nonumber
\end{eqnarray}
By Lemma \ref{t-derivative},
$$
\frac{\partial^{q}}{\partial
 t^{q}}u(t,x)=\int_{D}\int_{0}^{\infty}f(y)p_{s}^{D}(x,y)\frac{\partial^{q}}{\partial
 t^{q}}u^{\alpha /2}_{t}(s)dsdy.
$$
So we get
\begin{eqnarray}
\Delta ^{k}u(t,x)+(-1)^{k+1}\frac{\partial^{2m}
}{\partial t^{2m}}u(t,x)&=& \int_{D}\int_{0}^{\infty}f(y)
\left[(-1)^{k} \frac{\partial}{\partial s^{k}}u^{\alpha
/2}_{t}(s)\right.\nonumber\\
& +& \left. (-1)^{k+1}\frac{\partial^{2m}}{\partial
t^{2m}}u^{\alpha /2}_{t}(s)\right]
p_{s}^{D}(x,y)dsdy\nonumber\\
&=& 0,\ \mathrm{by \ Lemma\ \ref{density-pde}}.\nonumber
\end{eqnarray}

Now using dominated convergence theorem, we see that the boundary
conditions are satisfied.
\end{proof}

\textbf{Acknowledgments.} I would like to thank  Professor
Rodrigo Ba\~{n}uelos, my academic advisor, for  his guidance on
this paper.

\end{document}